\documentclass[12pt]{article}
\usepackage[latin1]{inputenc}
\usepackage{style}
\usepackage{geometry}
\usepackage{graphicx}
\usepackage{amsthm}
\usepackage[dvipsnames]{xcolor} 
\usepackage{enumitem}

 \makeatletter
\let\@fnsymbol\@alph
\makeatother

\theoremstyle{plain}
\newtheorem{theorem}{Theorem}[section]
\newtheorem{lemma}[theorem]{Lemma}
\newtheorem{notation}[theorem]{Notation}
\newtheorem{cor}[theorem]{Corollary}

\theoremstyle{definition}
\newtheorem{rem}[theorem]{Remark}

\numberwithin{equation}{section} 

\usepackage[onehalfspacing]{setspace} 

\newcommand{\frexp}[1]{\exp\bigl\{#1\bigr\}} 
\newcommand{\frphi}[1]{\Phi\bigl(#1\bigr)} 
\newcommand{\onebr}[1]{\one_{\{#1\}}} 

\newcommand{\sgn}{\mathrm{sgn}} 
\newcommand{\supnorm}[1]{\|#1\|_{\infty}}

\title{Stochastic Control of Drawdowns via Reinsurance under Random Inspection}
\author{
Kira Dudziak\footnote{University of Cologne, Department of Mathematics and Computer Science, Division of Mathematics, Cologne, Germany. E-mail: kira.dudziak@uni-koeln.de} \and Hanspeter Schmidli\footnote{Corresponding author. University of Cologne, Department of Mathematics and Computer Science, Division of Mathematics, Cologne, Germany. E-mail: hanspeter.schmidli@uni-koeln.de}
}

\begin{document}
\maketitle
\begin{abstract}
We consider a diffusion risk model where proportional reinsurance can be
bought. In order to stabilise the surplus process, one tries to keep the
drawdown, that is the difference of the surplus to its historical maximum, in
an interval $[0,d)$. The observation times of the drawdowns form a renewal
process. The retention levels can only be changed at the observation times
either. We show that an optimal strategy exists and how it is determined. We
illustrate the findings in the case of Poissonian observation times and
deterministic inter-observation times.
\end{abstract}

\medskip

\medskip

\noindent {Classification:} MSC: Primary 91B05; secondary 60G42, 93E20, 91G05;

\section{Introduction}

\subsection{Problem Formulation}

In order to appear reliable, insurance companies are interested in stabilising their surplus process and avoiding large losses. An object of interest is therefore the size of the drawdown, that is, the distance of the surplus to the last historical maximum.
In \cite{Brinker2022}, the time, in which the drawdown process of an insurer
exceeds a critical level $d > 0$, was considered and an optimal proportional
reinsurance strategy for minimising this time was determined. In reality, one
can only observe the process and adapt the strategy at discrete time
points. In this paper, we implement this idea by inspecting the process at the
arrival times of an ordinary renewal process. Our aim is to minimise
the number of time points
at which a critical drawdown level is observed. We state a dynamic programming
equation and show that the value function is the unique increasing and bounded
solution. Moreover, we calculate the distribution of the drawdown under
constant strategies explicitly.
The dynamic programming equation can be solved numerically, which we do
explicitly in the case of exponentially distributed and
deterministic interarrival times, respectively.

We work on a complete probability space $(\Omega,?F,!!P)$ containing all
the stochastic objects defined below. The surplus process
of an insurance portfolio is given by the diffusion approximation with
reinsurance
\[
X_{x_0}^B(t) = x_0 + \int_{0}^{t} \mu(B_s) \id s + \int_{0}^{t}
\sigma(B_s) \id W_s\;, \qquad t\geq 0\;,
\]
where $W = \{W_t\}_{t\geq 0}$ is a standard Brownian motion,
$\mu(b):=\eta-(1-b)\theta$, $\sigma(b):=\sigma b$ and $x_0$ is the surplus
level at the beginning of the observation period $t=0$. That is, we
consider proportional reinsurance. The parameters
$\eta$ and $\theta$ can be interpreted as the safety loading of the insurer
and reinsurer, respectively. In order that $B_t=0$ is not an optimal strategy,
we assume $\theta > \eta$, so reinsurance is more expensive than first insurance.
$B = \{B_t\}_{t\geq 0}$ is the reinsurance strategy
specified below. The surplus is not monitored continuously but is only
observed at the times $0 = T_0 < T_1 < T_2 < \cdots$. We model the observation
times by a renewal process $N=\{N_t\}_{t\geq 0}$ and denote the distribution
of the interarrival times by $F$ and its Laplace
transform by $\ell_{T}(r) = \E[\e^{-r T_1}]$. We assume that $N$ and $W$
are independent. The filtration $?F=\{?F_t\}_{t\geq 0}$ is the natural
filtration generated by $W$ and $N$, that is the smallest
right continuous filtration such that $W$ and $N$ are adapted. 

The insurer has the possibility to buy proportional reinsurance. But the
reinsurance treaty can only be changed at the observation times
$\{T_k\}$, too. That means
\[
B_t= \sum\limits_{k=0}^{\infty} b_{k} \one_{[T_{k},T_{k+1})}(t)= b_{N_t}\;,
\]
where $b_k$ are $?F_{T_k}$-measurable random variables with values in $[0,1]$
for all $k$. The set of all such strategies is
denoted by $?B$. For a chosen strategy $B \in ?B$, the \textit{(controlled)
running maximum} $M_{(m_0,x_0)}^B=\{M_{(m_0,x_0)}^B(t)\}_{t\geq 0}$ and the
\textit{(controlled) drawdown}
$\Delta_{(m_0,x_0)}^B=\{\Delta_{(m_0,x_0)}^B(t)\}_{t\geq 0}$ are given by 
\[
M_{(m_0,x_0)}^B(t)=\max\bigl\{m_0,\sup\limits_{0\leq s\leq
t}X_{x_0}^B(s)\bigr\}\;, \quad 
\Delta_{(m_0,x_0)}^B(t)= M_{(m_0,x_0)}^B(t)-X_{x_0}^B(t)\;.
\]
\begin{rem}
Note that the distribution of the drawdown process $\Delta_{(m_0,x_0)}^B$
depends on the difference $m_0 - x_0$ only. It therefore makes sense to write $\Delta_{m_0-x_0}^B$ instead. \end{rem}
Our aim is to stabilise the drawdown process in the what we call
\textit{non-critical area} $[0,d]$, that is, the drawdown size should not exceed
the level $d>0$. If we nevertheless enter the critical area $(d,\infty)$, what
cannot be prevented for all times, then we prefer to leave it before it is
observed at the inspection times $\{T_k\}$. We therefore consider
\[
v^B(z) = \E\Bigl[\sum\limits_{k=0}^{\infty}\e^{-r T_k}\one_{\{\Delta_{z}^B(T_k)>d\}}\Bigr]\;, \qquad z\geq 0\;,
\]
where we count the number of observations in the critical area and $z$ is the
initial drawdown. Note that we include the present observation.
We additionally include a preference factor $r>0$, such that the
present observation has a higher weight than an observation far in the
future. Since we want to minimise the number of observations, the
\textit{value function} is given by
\[
v(z) = \inf\limits_{B\in?B} v^B(z)\;, \qquad z\geq 0\;.
\]      
The paper is organised as follows. We start by giving some basic
properties of $v$, including the dynamic programming equation and uniqueness
of its solution. In Section~\ref{sec:DistDraw} we determine the distribution
of the drawdown at the next observation time. Two specific examples are
considered in Section~\ref{sec:scen}. Some technical details are given in the appendix.

\section{Dynamic Programming}

\subsection{First Results and Dynamic Programming}

We start with the following
\begin{lemma} \label{lemma:bound}
The value function $v$ is non-negative, increasing and bounded by $\lim_{z \to
\infty} v(z) = (1-\ell_{T}(r))^{-1}$. 
\end{lemma}
\begin{proof}
It is clear that $v(z)$ takes non-negative values only. For the upper bound we
observe
\[
v(z) \leq \E\Bigl[\sum\limits_{k=0}^{\infty}\e^{-r T_k}\Bigr] = \sum\limits_{k=0}^{\infty} [\ell_{T}(r)]^k = \frac{1}{1-\ell_{T}(r)}\;.
\]
Moreover, note that $z\mapsto \Delta_z^B(t)$ is (pathwise) increasing for
every $t\geq 0$ and $B\in?B$. It is then easy to see that for $y\leq z$ we get
$v(y)\leq v^B(y)\leq v^B(z)$ for every $B\in?B$. By taking the infimum over
$B\in?B$ on the right hand side, we conclude that $v$ is increasing.
The limit as $z \to \infty$ follows from Bellman's equation below, see
Corollary~\ref{cor:limit}.
\end{proof}
We can look at the observation times $\{T_k\}$ as regeneration
times. Thus, for us the process up to time $T_1$ is of interest in which the
chosen strategy is constant. We therefore consider now processes with a
constant strategy.
\begin{notation}
For $b\in[0,1]$ we write $X_{x_0}^b$,
$M_{(m_0,x_0)}^b$ and $\Delta_{x_0}^b$ for
the constant strategy $B_t = b$. 
\end{notation}


The following result is crucial for the characterisation of the value
function.
\begin{theorem}\label{thm:dynamicprogramming}
The value function $v$ fulfils the dynamic programming equation
\begin{equation}\label{dynamicprogramming}
v(z) = \one_{\{z>d\}}+ \inf_{b\in[0,1]}  \E[\e^{-r T_1}v(\Delta_{z}^b(T_1))]\;, \quad z\geq 0\;. 
\end{equation}
\end{theorem}
\begin{proof}
Let $B\in?B$ be an admissible strategy with $B_t=b$ for $t<T_1$ and let $z\geq
0$. We observe
\begin{eqnarray*}
v^B(z) &=& \onebr{z>d} + \E\Bigl[\sum_{k=1}^{\infty} \e^{-rT_k}\onebr{\Delta_z^B(T_k)>d}\Bigr] \\
&=& \onebr{z>d} + \E\Bigl[\e^{-rT_1} \E\Bigl[\sum\limits_{k=1}^{\infty} \e^{-r(T_k-T_1)}\onebr{\Delta_z^B(T_k)>d}\bigm| ?F_{T_1}\Bigr]\Bigr] \\
&=& \onebr{z>d} + \E\Bigl[\e^{-rT_1} \E\Bigl[\sum\limits_{k=0}^{\infty} \e^{-r\tilde{T}_k}\onebr{\Delta_{\Delta_z^{b}(T_1)}^{\tilde{B}}(\tilde{T}_k)>d}\bigm|  T_1 \Bigr]\Bigr] \\
&=& \onebr{z>d} + \E[\e^{-rT_1}v^{\tilde{B}}(\Delta_z^{b}(T_1))] \geq
\onebr{z>d} + \inf\limits_{b\in[0,1]}\E[\e^{-rT_1}v(\Delta_z^{b}(T_1))]\;, 
\end{eqnarray*} 
where $\tilde{T_k} :=T_{k+1}-T_1$ for $k\in!!N_0$ and $\tilde{B}_t:=B_{t+T_1}$ for $t\geq 0$. 
Taking the infimum over $B\in?B$ on the left hand side yields
\begin{equation}
v(z) \geq \onebr{z>d} + \inf\limits_{b\in[0,1]}\E[\e^{-rT_1}v(\Delta_z^{b}(T_1))]\;. \label{eq:firstineq}
\end{equation}
Choose $\epsilon > 0$ and $\delta > 0$.  
Since the image of $v$ is bounded, there are only finitely many jumps larger
than $\delta$. Thus, we can choose $m \in !!N$ and $0 = \xi_0 < \xi_1 < \cdots < \xi_m$, such
that $\{z: v(z+) - v(z-) \ge \delta\} \subset \{\xi_0,\xi_1, \ldots, \xi_m\}$. Moreover, for every $l\in\{0,\ldots,m\}$ we can choose $\{z_k^l\}_{k\in!!N_0}$, such that $\xi_l=z_0^l < z_1^l < \cdots < \xi_{l+1}$, $\lim_{k\to\infty} z_k^l = \xi_{l+1}$ and $v(z_k^l)-v(z) < \delta$ for $z_{k-1}^l < z \leq z_k^l$. We define 
\[
k(z) = \begin{cases}
z\,, & z\in\{\xi_0,\ldots,\xi_m\}\,, \\
z_k^l\,, & z\in (z_{k-1}^l,z_k^l]\,.
\end{cases}
\]
For every $l\in\{0,\ldots,m\}$ we can choose a strategy $B_l^{\epsilon}$, such that $v^{B_l^\epsilon}(\xi_l) < v(\xi_l) + \epsilon$. Moreover, for every $k\in !!N_0$ we can choose $B_l^{\epsilon,k}$, such that $v^{B_l^{\epsilon,k}}(z_k^l) < v(z_k^l) + \epsilon$. We define
\[
B^{\epsilon}(z) = \begin{cases}
B_l^{\epsilon}\,, & k(z)=\xi_l\,, \\
B_l^{\epsilon,k}\,, & k(z)= z_k^l\,.
\end{cases}
\]
Now let $b\in[0,1]$ be arbitrary and consider the strategy $B_t = b\one_{\{t<T_1\}} + B^{\epsilon}_{t-T_1}(\Delta_z^b(T_1))\one_{\{t\geq T_1\}}$.
Analogously to above we find 
\begin{eqnarray*}
v(z) &\leq& v^B(z) = \onebr{z>d} +
\E[\e^{-rT_1}v^{B^{\epsilon}(\Delta_z^b(T_1))}(\Delta_z^b(T_1))] \le \onebr{z>d} +
\E[\e^{-rT_1}v^{B^{\epsilon}(\Delta_z^b(T_1))}(k(\Delta^b_z(T_1)))]
\\
&\leq& \onebr{z>d} + \E[\e^{-rT_1} v(k(\Delta^b_z(T_1)))]+\epsilon
\le \onebr{z>d} + \E[\e^{-rT_1} v (\Delta_z^b(T_1))] + \delta + \epsilon\;. 
\end{eqnarray*}
By letting $\delta \to 0$, $\epsilon\to 0$ and by taking the infimum over $b\in[0,1]$, we conclude
\begin{equation*}
v(z) \leq \onebr{z>d} + \inf\limits_{b\in[0,1]}\E[\e^{-rT_1}v(\Delta_z^{b}(T_1))]\;.
\end{equation*}
Together with \eqref{eq:firstineq}, this proves the assertion.
\end{proof}
\begin{lemma} \label{lemma:distribution}
For a fixed retention level $b$, the distribution of the drawdown at time $t$
is given by the distribution function
\begin{equation}
F_{\Delta_z^b(t)}(\delta) = !!P[\Delta_z^b(t)\leq \delta] = \frphi{\frac{\delta-z+t\mu(b)}{\sqrt{t}\sigma(b)}}-\frexp{-2\delta\frac{\mu(b)}{\sigma^2(b)}}\frphi{\frac{-\delta-z+t\mu(b)}{\sqrt{t}\sigma(b)}}\;. \label{drawdowndistribution}
\end{equation}
The density is of the form
\begin{eqnarray}
f_{\Delta_z^b(t)}(\delta) &=& \Bigl[\frac1{\sqrt{2\pi t} \sigma(b)} \Bigl(
\frexp{-\frac{(\delta + z+ t \mu(b))^2 - 4 t \mu(b) z}{2 t
\sigma^2(b)}}+ \frexp{-\frac{(\delta-z+ \mu(b)t)^2}{2 t \sigma^2(b)}}\Bigr)
\nonumber\\
&& {}+ \frac{2 \mu(b)}{\sigma^2(b)}
\Phi\Bigl(\frac{t \mu(b) - \delta-z}{\sqrt{t}\sigma(b)}\Bigr)
\e^{-2\delta \mu(b)/\sigma^2(b)}\Bigr] \one_{\delta > 0}\;. \label{density}
\end{eqnarray}
\end{lemma}
We postpone the proof to Section \ref{sec:DistDraw}. 
\begin{cor}\label{cor:limit}
We have $\lim_{z \to \infty} v(z) = [1-\ell_{T}(r)]^{-1}$.
\end{cor}
\begin{proof}
Let $\zeta = \lim_{z \to \infty} v(z)$ and fix $\epsilon > 0$ and $\delta >
0$. Then there is 
$z_0$, such that $v(z) > \zeta - \epsilon$ for any $z \ge z_0.$ 
Next, we show that there is $z_1>z_0$, such that $\Prob[\Delta_z^b(T_1) \le z_0] < \delta$ for
$z \ge z_1$ and any $b \in [0,1]$. Note that we can find $a>0$ with $\Prob[T_1>a]<\delta/2$. Choose 
\[
z_1 > z_0 + a\eta - \sqrt{a}\sigma\Phi^{-1}\bigl(\frac{\delta}{2F_{T_1}(a)}\bigr)\;.
\]
We can assume that $\delta/(2F_{T_1}(a)) < 1/2$, i.e.~$z_1>z_0+a\eta$. By using \eqref{drawdowndistribution}, this yields for $z\geq z_1$, $b\in[0,1]$ and $t\in[0,a]$:
\begin{eqnarray*}
\Prob[\Delta_z^b(t) \leq z_0] &\leq& \Phi\bigl(\frac{z_0+t\mu(b)-z}{\sqrt{t}\sigma(b)}\bigr) \leq \Phi\bigl(\frac{z_0+a\eta - z}{\sqrt{a}\sigma}\bigr) \\
&\leq& \Phi\bigl( \frac{z_0+a\eta-[z_0 + a\eta - \sqrt{a}\sigma\Phi^{-1}(\delta/(2F_{T_1}(a)))]}{\sqrt{a}\sigma}\bigr) \\
&=& \frac{\delta}{2F_{T_1}(a)}\;.
\end{eqnarray*}
We conclude
\begin{eqnarray*}
\Prob[\Delta_z^b(T_1) \leq z_0 ] &=& \Prob[\Delta_z^b(T_1) \leq z_0, T_1>a] + \Prob[\Delta_z^b(T_1) \leq z_0, T_1\leq a] \\
&<& \frac{\delta}{2} + \int_0^{a} \Prob[\Delta_z^b(t) \leq z_0] \id F_{T_1}(t) \\
&\leq& \frac{\delta}{2} + \int_0^{a} \frac{\delta}{2F_{T_1}(a)} \id F_{T_1}(t) \\
&=& \delta\;.
\end{eqnarray*}
Thus, for $z \ge z_1\vee d$
\[ v(z) \ge 1 + \E[\e^{- r T_1}] (\zeta - \epsilon)(1-\delta) \;.\]
Letting $z \to \infty$, $\zeta \ge 1 + \E[\e^{- r T_1}] (\zeta -
\epsilon)(1-\delta)$ and because $\delta$ and $\epsilon$ are arbitrary, $\zeta
\ge  1 + \ell_{T}(r) \zeta$. Thus, $\zeta \ge [1-\ell_{T}(r)]^{-1}$.
\end{proof}
Next, we show that the value function $v$ is the unique solution to the
dynamic programming equation \eqref{dynamicprogramming}.

\subsection{Existence and Uniqueness of a Solution}

Denote by $?S$ the space of all bounded functions $f$ with values in $[0,[1-\ell_{T}(r)]^{-1}]$. Then, $(?S,\supnorm{\cdot})$ is a Banach space.
\begin{theorem}
The value function $v$ is the unique solution $v\in?S$ to the dynamic programming equation \eqref{dynamicprogramming}.
\end{theorem}
\begin{proof}
We define the operator $?T^b$ by
\begin{equation}
?T^b(f)(z):=\onebr{z>d}+\E[\e^{-r T_1}f(\Delta_z^b(T_1))], \qquad f\in?S, \quad z\geq 0 \label{operator}
\end{equation}
and let $?T(f)(z) = \inf_{b \in [0,1]} ?T^b(f)(z)$. We first note that
\[ |?T(f)(z)| \le 1 + \E[\e^{-r T_1} |f(\Delta_z^b(T_1))|] \le 1 +
\ell_{T}(r)/(1 - \ell_{T}(r)) = (1 - \ell_{T}(r))^{-1}\;,\]
and therefore $?T(f) \in ?S$ for $f \in ?S$.
Let $f,g\in?S$, $z\geq 0$ and $b \ge 0$. We have
\begin{eqnarray*}
?T(f)(z)-?T^b(g)(z) &\le& ?T^b(f)(z)-?T^b(g)(z) = \E[\e^{-r  T_1}  (f(\Delta_z^b(T_1))-g(\Delta_z^b(T_1)))]\\
&=&  \int_0^\infty \e^{-rt}\int_0^\infty (f(y)-g(y)) \id F_{\Delta_z^b(t)}(y) \id F(t) \\
&\leq& \int_0^\infty \e^{-rt}\int_0^\infty \abs{f(y)-g(y)}\id F_{\Delta_z^b(t)}(y) \id F(t) \leq \ell_{T}(r) \supnorm{f-g}\;.
\end{eqnarray*}
Maximising over $b$ gives $?T(f)(z)-?T(g)(z) \le \ell_{T}(r)
\supnorm{f-g}$. Analogously, $?T(g)(z)-?T(f)(z) \le \ell_{T}(r)
\supnorm{f-g}$, yielding $|?T(f)(z)-?T(g)(z)| \le \ell_{T}(r)
\supnorm{f-g}$.
By taking the supremum over $z\in[0,\infty)$ on the left hand side, we
conclude $\supnorm{?T^b(f)-?T^b(g)}\leq \ell_{T}(r) \supnorm{f-g}$.
This proves that $?T$ is a contraction. By Banach's fixed point theorem, there
exists a unique solution $v\in?S$ to \eqref{dynamicprogramming} 
which, according to Theorem \ref{thm:dynamicprogramming}, has to be the value
function.
\end{proof}
\begin{rem} \label{rem:Banach}
Banach's fixed point theorem also provides an approach for determining the value
function numerically. For $f_0\in?S$ define $f_{n+1}=?T(f_n)$ for all
$n\in!!N$. This series converges to the fixpoint exponentially fast, which is the value function. We will exploit this in order to calculate $v$ numerically.
\end{rem} 

\begin{theorem}
There exists a measurable function $b^{*}:[0,\infty)\to[0,1]$ such that for all $z\geq 0$
\begin{equation}
\inf_{b\in[0,1]} \E[\e^{-rT_1}v(\Delta_z^b(T_1))] = \E[\e^{-rT_1}v(\Delta_z^{b^{*}(z)}(T_1))]\;. \label{eq:minimiser}
\end{equation}
\end{theorem}
\begin{proof}

We claim that the continuous function 
\[ [0,\infty) \times (0,1] \to !!R\;,\quad (z,b) \mapsto \E[\e^{-r T_1}
v(\Delta_z^b(T_1))] = \E\bigl[\e^{-r T_1}\int_0^\infty v(\delta)
f_{\Delta_z^b(T_1)}(\delta) \id \delta\bigr]\]
can be extended to a continuous
function $a:[0,\infty) \times [0,1] \to !!R$ by
\[
a(z,b) := \begin{cases}
\E[\e^{-rT_1} v(\Delta_z^b(T_1))]\;, & b>0\;, \\
\lim_{b\downarrow 0} \E[\e^{-rT_1} v(\Delta_z^b(T_1))]\;, & b=0\;.
\end{cases}
\]

For $t=T_1$ we consider the three terms in \eqref{density} $f_{\Delta_z^b(T_1)}(\delta) =
f_1(\delta) + f_2(\delta) + f_3(\delta)$ separately. That is, we calculate the
limits of $\int_0^\infty f_k(\delta) v(\delta) \id \delta$ as $b
\downarrow 0$. In the following we assume that $\mu(b) < 0$.
We start with $f_2$. The expression can be written as
\[ \E[v(\sigma(b) \sqrt{T_1} ?Z + z - T_1 \mu(b)) \one_{\sigma(b) \sqrt{T_1}
?Z > T_1 \mu(b)-z}]\;,\]
where $?Z$ is a standard normally distributed random variable. Splitting the
expression in $\{?Z \ge \theta \sqrt{T_1}/?\sigma\}$ and $\{?Z < \theta
\sqrt{T_1}/?\sigma\}$, we get that the limit $b \downarrow 0$ is
\begin{equation}
v((z +\{\theta - \eta\} T_1) +) \Prob[?Z > \sqrt{T_1}\theta/\sigma] +
v((z+\{\theta - \eta\} T_1) -) \Prob[?Z < \sqrt{T_1}\theta/\sigma]\;. \label{eq:limitbtozero}
\end{equation}
Next, consider $f_3$. Note that $-(\sqrt{t} \sigma(b))^{-1}(z + \delta - t
\mu(b))$ tends to $-\infty$. Using $\Phi(x) \sim - \e^{-x^2/2}/(x \sqrt{2
\pi})$ as $x \to -\infty$, yields
\begin{eqnarray*}
f_3(\delta) &\approx& \frac{2 \mu(b)t}{\sqrt{2 \pi t}\sigma(b)(z + \delta - \mu(b)
t)} \exp\bigl\{-\frac{(\delta + z - t \mu(b))^2 + 4 \delta t\mu(b))}{2 t
\sigma^2(b)}\bigr\} \one_{\delta > 0}\\
&=& \frac{2 \mu(b)t}{\sqrt{2 \pi t}\sigma(b)(z + \delta - \mu(b)
t)} \exp\bigl\{-\frac{(\delta + z + t \mu(b))^2 - 4 t\mu(b)z}{2 t
\sigma^2(b)}\bigr\} \one_{\delta > 0}\;.
\end{eqnarray*}
Thus, the terms connected to $f_1 + f_3$ can be written as
\begin{eqnarray*}
\lefteqn{\E\bigl[\bigl(1 + \frac{2 \mu(b) T_1}{\sqrt{T_1} \sigma(b) ?Z - 2 \mu(b)
T_1 } \bigr) v( \sqrt{T_1} \sigma(b) ?Z - \mu(b) T_1-z) \one_{\sqrt{T_1}
\sigma(b) ?Z > \mu(b) T_1+z} \bigr] \exp\bigl\{\frac{2 \mu(b) z}{\sigma^2(b)}
\bigr\}}\hskip1cm\\
&=& \E\bigl[\frac{\sqrt{T_1} \sigma(b) ?Z}{\sqrt{T_1} \sigma(b) ?Z - 2 \mu(b)
T_1 } v( \sqrt{T_1} \sigma(b) ?Z - \mu(b) T_1-z) \one_{\sqrt{T_1}
\sigma(b) ?Z > \mu(b) T_1+z} \bigr] \exp\bigl\{\frac{2 \mu(b) z}{\sigma^2(b)}
\bigr\}\;.
\end{eqnarray*}
Note that $\sqrt{T_1} \sigma(b) ?Z - 2 \mu(b) T_1 \ge z - \mu(b) T_1 > 0$ in
the integration area. It now follows for $b\downarrow 0$ that the expression considered tends to
zero if $z > 0$ and also if $z = 0$. This shows that $a(z,0)$ exists.

Now there exists a minimiser $b_a(z)$ such that
$a(z,b_a(z))=\inf_{b\in[0,1]}a(z,b)$. By \cite{wagner}, we can choose the
function $b_a(z)$ measurably.
Note that $a(z,0) \geq \E[\e^{-r T_1} v((z+(\theta - \eta)T_1)-)]$ by \eqref{eq:limitbtozero} and the observations above. This yields
\[
\inf_{b \in [0,1]} \E[\e^{-r T_1} v(\Delta_z^b(T_1))] = \min \{ a(z,b_a(z)), \E[\e^{-r T_1} v(\Delta_z^0(T_1))] \}\;,
\]
if $v$ is left-continuous. In this case, the minimiser
\begin{equation}
b^{*}(z) = b_a(z) \one_{a(z, b_a(z)) \leq \E[\e^{-r T_1} v(\Delta_z^0(T_1))]} \label{minimiser}
\end{equation}
fulfils \eqref{eq:minimiser}. Suppose that $v$ is not left-continuous. Choose $z$, such that $v(z) - v(z-)$ is maximal. This is possible because $v$ is bounded. Due to the continuity of $(z,b) \mapsto \E[\e^{-r T_1} v(\Delta_z^b(T_1))]$ for $b>0$ and the dynamic programming equation, $v(z)>v(z-)$ can only be true if
\[
v(z) = \one_{z>d} + \min \{ a(z,0), \E[\e^{-r T_1} v(\Delta_z^0(T_1))] \}
\] 
and
\begin{eqnarray*}
v(z-) &=& \lim_{\epsilon\downarrow 0} \bigl\{\one_{z-\epsilon>d} + \min \{ a(z-\epsilon,0), \E[\e^{-r T_1} v(\Delta_{z-\epsilon}^0(T_1))] \} \bigr\} \\
&=& \one_{z>d} + \min \{ a(z,0), \E[\e^{-rT_1} v(\Delta_z^0(T_1)-)]\}
= \one_{z>d} + \E[\e^{-rT_1} v(\Delta_z^0(T_1)-)]\;.
\end{eqnarray*}
But this implies that
\begin{eqnarray*}
0 \le v(z)-v(z-) &=& \min \{ a(z,0), \E[\e^{-r T_1} v(\Delta_z^0(T_1))] \} - \E[\e^{-rT_1} v(\Delta_z^0(T_1)-)] \\
&\leq& \E[\e^{-rT_1}(v(\Delta_z^0(T_1))-v(\Delta_z^0(T_1)-))]\;.
\end{eqnarray*}
Since $!!P[v(\Delta_z^0(T_1))=v(\Delta_z^0(T_1)-)] = 1$ would violate $v(z) > v(z-)$, we find
\[
v(z)-v(z-) \leq \E[\e^{-rT_1}(v(\Delta_z^0(T_1))-v(\Delta_z^0(T_1)-))] < \E[v(\Delta_z^0(T_1))-v(\Delta_z^0(T_1)-)]\;,
\]
which contradicts that $v(z)-v(z-)$ was maximal. We conclude that $v(z)$ is left-continuous. This shows that $b^{*}(z)$ defined in \eqref{minimiser} is a minimiser.
\end{proof}

\begin{theorem}
The strategy $B^{*}$ given by
\[
B_t^{*} = b^{*}(\Delta_{z}^{B^{*}}(T_{N_t}))\;,
\]
is optimal, in the sense that $v(z)=\inf_{B\in?B}v^{B}(z)=v^{B^{*}}(z)$ for all $z\geq 0$.
\end{theorem}
\begin{proof}
Consider the process 
\[ V_n = \sum_{k=0}^{n-1} \e^{-r T_k}\onebr{\Delta_z^{B^{*}}(T_k)>d} + \e^{-r
T_n} v(\Delta_z^{B^{*}}(T_n))\;.\]
Then,
\begin{eqnarray*}
\E[V_{n+1} \mid ?F_{T_n}] &=&  \sum_{k=0}^{n-1} \e^{-r
T_k} \onebr{\Delta_z^{B^{*}}(T_k)>d} + \e^{-r T_n}\bigl\{
\onebr{\Delta_z^{B^{*}}(T_n)>d} + \E[\e^{-r (T_{n+1}-T_n)}
v(\Delta_z^{B^{*}}(T_{n+1}))\mid ?F_{T_n}]\bigr\}\\ 
&=& \sum_{k=0}^{n-1} \e^{-r
T_k} \onebr{\Delta_z^{B^{*}}(T_k)>d} + \e^{-r T_n} v(\Delta_z^{B^{*}}(T_n)) =
V_n\;, 
\end{eqnarray*}
where we used \eqref{dynamicprogramming} in the last but one equality. Thus,
$V_n$ is a martingale. This gives
\begin{eqnarray*}
v(z) &=& V_0 = \lim_{n \to \infty} \E[V_n] = \lim_{n \to \infty} \E\Bigl[\sum_{k=0}^{n-1} \e^{-r T_k}\onebr{\Delta_z^{B^{*}}(T_k)>d} + \e^{-r
T_n} v(\Delta_z^{B^{*}}(T_n))\Bigr] \\
&=& \E\Bigl[ \sum_{k=0}^\infty
\e^{-r T_k}\onebr{\Delta_z^{B^{*}}(T_k)>d}\Bigr] + 0\;,
\end{eqnarray*}
where we used monotone convergence and bounded convergence, respectively. This
shows that $v(z)$ is the value of the proposed strategy. 
\end{proof}

\section{Distribution of the Drawdown}\label{sec:DistDraw}

The dynamic programming equation \eqref{dynamicprogramming} will be essential
for the numerical calculation of the value function and the optimal strategy. We would therefore like to rewrite the equation to
\begin{equation}
v(z) = \onebr{z>d} + \inf_{b\in[0,1]} \int_0^\infty \e^{-rt}\! \int_0^\infty v(\delta) \id F_{\Delta_z^b(t)}(\delta) \id G(t)\;,
\end{equation}
where $G$ denotes the distribution function of $T_1$. Thus, we need to determine the distribution of $\Delta_z^b(t)$. If $b=0$, we immediately see that the drawdown process is deterministic and $!!P[\Delta_z^0(t)\leq \delta] = \one_{t\leq(\delta-z)/(\theta-\eta)}$. In the following, we calculate the distribution of $\Delta_z^b(t)$ for $b\in(0,1]$. We start with the case where the initial drawdown is zero.
\begin{lemma}
The density of $(\Delta_0^b(t),X_0^b(t))$ is
\begin{equation}
f_{(\Delta_0^b(t),X_{0}^b(t))}(\delta,x)= \frac{2(2\delta +x)}{\sqrt{2\pi t^3}
\sigma(b)^3}\frexp{-\frac{(2\delta+x)^2}{2t\sigma^2(b)}+\frac{\mu(b)}{\sigma^2(b)}x-\frac{t\mu(b)^2}{2\sigma^2(b)}}\one_{\{x+\delta \geq 0\}}\one_{\{\delta\geq 0\}}\;.
\end{equation}
\end{lemma}
\begin{proof}
The joint density of a Brownian motion with drift $(X_t)_{t\geq 0}=(\mu t+W_t)_{t\geq 0}$ and its running maximum $(M_t)_{t\geq 0}=(\sup_{s\in[0,t]}X_s)_{t\geq 0}$ is given by
\[
f_{(M,X)}(m,x)= \frac{2(2m-x)}{\sqrt{2\pi t^3}}\frexp{-\frac{(2m-x)^2}{2t}+\mu x-\frac{\mu^2t}{2}}\one_{\{m\geq 0\}}\one_{\{x\leq m\}}\;,
\]
see for example \cite[Equation (13.10)]{Rogers_Williams_2000}. This
yields for $(M_{(0,0)}^b(t),X_{0}^b(t))$
\[
f_{(M_{(0,0)}^b(t),X_{0}^b(t))}(m,x)= \frac{2(2m-x)}{\sqrt{2\pi t^3}\sigma(b)^3}\frexp{-\frac{(2m-x)^2}{2t\sigma^2(b)}+\frac{\mu(b)}{\sigma^2(b)}x-\frac{t\mu(b)^2}{2\sigma^2(b)}}\one_{\{m\geq 0\}}\one_{\{x\leq m\}}\;.
\]
By change of variable we get for the joint density of
$(\Delta_0^b(t),X_0^b(t)) = (M_{(0,0)}^b(t)-X_0^b(t),X_0^b(t))$
\[
f_{(\Delta_0^b(t),X_{0}^b(t))}(\delta,x)= \frac{2(2\delta +x)}{\sqrt{2\pi t^3}
\sigma(b)^3}\frexp{-\frac{(2\delta+x)^2}{2t\sigma^2(b)}+\frac{\mu(b)}{\sigma^2(b)}x-\frac{t\mu(b)^2}{2\sigma^2(b)}}\one_{\{x+\delta \geq 0\}}\one_{\{0\leq \delta\}}\;,
\]
which is the assertion.
\end{proof}
\begin{proof}[Proof of Lemma \ref{lemma:distribution}]
The initial drawdown is given by $\Delta_z^b(0)=z=m_0-x_0>0$
and let $x_0 = X_0^b(0) = 0$. If the process reaches the maximum in $(0,t]$,
that is $X_0^b(t) + \Delta_0^b(t) \ge z$, then $\sup_{s\in[0,t]}X_0^b(t)-X_0^b(t)$ is the
drawdown in $t$. The corresponding density on $\{M_0^b(t) \ge z\}$ is
then obtained by integrating the joint density over
$(z-\delta,\infty)$ with respect to $x$
\begin{eqnarray*}
f_{\Delta_z^b(t) ; m\ge z}(\delta) &=& \bigl[\frac2{\sqrt{2 \pi t} \sigma(b)}
\frexp{-\frac{(\delta + z- t \mu(b))^2 + 4 t \mu(b) \delta}{2 t \sigma^2(b)}}\\
&&\hskip2cm {}+ \frac{2 \mu(b)}{\sigma(b)^2} \Phi\bigl(\frac{t \mu(b) -
\delta-z}{\sqrt{t}\sigma(b)}\bigr) \e^{-2\delta \mu(b)/\sigma^2(b)}\bigr]
\one_{\delta \ge 0}\\
&=& \bigl[\frac2{\sqrt{2 \pi t} \sigma(b)}
\frexp{-\frac{(\delta + z+ t \mu(b))^2 - 4 t \mu(b) z}{2 t \sigma^2(b)}}\\
&&\hskip2cm {}+ \frac{2 \mu(b)}{\sigma(b)^2} \Phi\bigl(\frac{t \mu(b) -
\delta-z}{\sqrt{t}\sigma(b)}\bigr) \e^{-2\delta \mu(b)/\sigma^2(b)}\bigr]
\one_{\delta \ge 0}
\;.
\end{eqnarray*}
If the maximum in $(0,t]$ is not reached, then the drawdown is $z - X_0^b(t)$.
Integrating over $(0,z-x]$ with respect to $\delta$ yields the
density of $X_t^b(0)$ on $\{M_0^b(t) < z\}$
\[ f_{X_0^b(t) ; m < z}(x) = \frac1{\sqrt{2\pi t} \sigma(b)}
\bigl(\exp\bigl\{-\frac{x^2}{2 t \sigma^2(b)}\bigr\} - \exp\bigl\{-
\frac{(2z-x)^2}{2 t \sigma^2(b)}\bigr\}  \bigr)
\frexp{\frac{\mu(b) x}{\sigma^2(b)}-\frac{t\mu(b)^2}{2\sigma^2(b)}}\one_{x < z}\;.\]
This gives the density for the drawdown $\Delta_z^b(t) = z - X_0^b(t)$
on $\{M_0^b(t) < z\}$
\begin{eqnarray*}
f_{\Delta_z^b(t) ; m < z}(\delta) &=& \frac1{\sqrt{2\pi t} \sigma(b)}
\bigl(\exp\bigl\{-\frac{(z-\delta)^2}{2 t \sigma^2(b)}\bigr\} - \exp\bigl\{-
\frac{(z+\delta)^2}{2 t \sigma^2(b)}\bigr\}  \bigr)
\frexp{\frac{\mu(b) (z-\delta)}{\sigma^2(b)}-\frac{t\mu(b)^2}{2\sigma^2(b)}}
\one_{\delta > 0}\\
&=& \frac1{\sqrt{2\pi t} \sigma(b)}
\bigl(\frexp{-\frac{(z-\delta-\mu(b)t)^2}{2 t \sigma^2(b)}}
- \frexp{-\frac{(z+\delta+\mu(b)t)^2}{2 t \sigma^2(b)}+\frac{2 \mu(b) z}{\sigma^2(b)}}\bigr)\one_{\delta > 0}\;.
\end{eqnarray*}
Adding the two densities yields the density \eqref{density} of
$\Delta_z^b(t)$.
By integrating with respect to $\delta$ we find the distribution function $F_{\Delta_z^b(t)}$ of $\Delta_z^b(t)$:
\begin{eqnarray*}
F_{\Delta_z^b(t)}(\delta) = !!P[\Delta_z^b(t)\leq \delta] = \frphi{\frac{\delta-z+t\mu(b)}{\sqrt{t}\sigma(b)}}-\frexp{-2\delta\frac{\mu(b)}{\sigma^2(b)}}\frphi{\frac{-\delta-z+t\mu(b)}{\sqrt{t}\sigma(b)}}\;.
\end{eqnarray*}
\end{proof}
\section{Two Different Scenarios for the Renewal Process}\label{sec:scen}

\subsection{Poissonian Observations}

In this section, we assume that the renewal process $N$ is given by a Poisson
process with Poisson parameter $\rho$. The inspection of the process therefore
occurs at random time points. 

\subsubsection{Rewriting the Dynamic Programming Equation}

First, we note that, according to Lemma \ref{lemma:bound}, the value function is bounded from above by $(\rho+r)/r$. 
We rewrite the dynamic programming equation \eqref{dynamicprogramming} into an equation of the form
\begin{equation}
v(z) = \onebr{z>d} + \inf_{b\in[0,1]} \bigl\{ \int_0^\infty v(\delta)  w^b(\delta,z) \id \delta \bigr\} \label{dynprorewritten}
\end{equation}
and develop an algorithm in order to calculate $v(z)$ in a numerical way. If $b=0$, we have
\begin{eqnarray}
\E[\e^{-r T_1}v(\Delta_z^0(T_1))] &=& \int_0^\infty \rho\e^{-(r+\rho)t}v(z+(\theta-\eta)t) \id t =\int_z^\infty \frac{\rho}{\theta-\eta}\e^{-(r+\rho)(\delta-z)/(\theta-\eta)}v(\delta) \id \delta \nonumber\\
&=& \int_0^\infty w^0(\delta,z)  v(\delta) \id \delta\;,
\end{eqnarray}
giving
\[
w^0(\delta,z):= \frac{\rho}{\theta-\eta}\frexp{-(\delta-z) \frac{r+\rho}{\theta-\eta}}\one_{\{\delta> z\}}\;.
\]
We need the following result for the calculations in the case $b>0$.
\begin{lemma} \label{lemma:integralsinverseGaussian}
Let $s>0$, $\alpha,\mu\in!!R$ and $\sigma>0$. We have the following identities:
\begin{eqnarray*}
&&\int_0^\infty \e^{-st} \frac{1}{\sqrt{2\pi
t^3}\sigma}\frexp{-\frac{(\alpha-\mu t)^2}{2t\sigma^2}} \id t =
\abs{\alpha}^{-1} \e^{q}\;, \\ 
&&\int_0^\infty \e^{-st} \frac{1}{\sqrt{2\pi
t}\sigma}\frexp{-\frac{(\alpha-\mu t)^2}{2t\sigma^2}} \id t =
\frac{1}{\sqrt{\mu^2+2s\sigma^2}} \e^{q} \qquad\text{and} \\
&&\int_0^\infty \e^{-st} \frphi{\frac{\alpha-\mu t}{\sigma\sqrt{t}}} \id t = \frac{1}{s} (\one_{\{\alpha>0\}}+\frac{1}{2}\one_{\{\alpha=0\}})-\frac{1}{2s} \bigl(\sgn(\alpha)+\frac{\mu}{\sqrt{\mu^2+2s\sigma^2}}\bigr) \e^{q}\;,\\ 
\end{eqnarray*}
where
\[
q=\frac{\alpha \mu-\abs{\alpha} \sqrt{\mu^2+2s\sigma^2}}{\sigma^2}\;.
\]
\end{lemma}
The proof is given in the appendix.

Now we note that for $b\in(0,1]$ and $t\geq 0$, there exists a density function $f_{\Delta_z^b(t)}$ given by \eqref{density} for the distribution of $\Delta_z^b(t)$. We can therefore write
\begin{eqnarray*}
\E[\e^{-r T_1}v(\Delta_z^b(T_1))] &=& \int_0^\infty \rho\e^{-(r+\rho)t}\int_0^\infty v(\delta)  f_{\Delta_z^b(t)}(\delta) \id \delta \id t \\
&=& \int_0^\infty \rho v(\delta)\int_0^\infty \e^{-(r+\rho)t}f_{\Delta_z^b(t)}(\delta) \id t \id \delta\;.
\end{eqnarray*}
Let $\zeta(b):=\sqrt{\mu(b)^2+2(r+\rho)\sigma^2(b)}$. By using Lemma \ref{lemma:integralsinverseGaussian} we find
\begin{eqnarray*}
\rho^{-1}  w^b(\delta,z)&:=&\int_0^\infty \e^{-(r+\rho)t}f_{\Delta_z^b(t)}(\delta) \id t \\
&=& \int_0^\infty \e^{-(r+\rho)t} \frac{1}{\sqrt{2\pi t}\sigma(b)}\frexp{-\frac{(z-\delta-t\mu(b))^2}{2t\sigma^2(b)}} \id t \\
&&\quad {}+ \frac{2\mu(b)}{\sigma^2(b)}\frexp{-2\delta\frac{\mu(b)}{\sigma^2(b)}} \int_0^\infty \e^{-(r+\rho)t} \frphi{\frac{-\delta-z+t\mu(b)}{\sqrt{t}\sigma(b)}} \id t \\
&&\quad {}+ \frexp{-2\delta\frac{\mu(b)}{\sigma^2(b)}}  \int_0^\infty \e^{-(r+\rho)t} \frac{1}{\sqrt{2\pi t}\sigma(b)}\frexp{-\frac{(\delta+z-t\mu(b))^2}{2t\sigma^2(b)}} \\
&=& \zeta(b)^{-1} \frexp{\frac{(z-\delta)\mu(b)-\abs{z-\delta} \zeta(b)}{\sigma^2(b)}} \\ 
&\quad&+\frac{2\mu(b)}{\sigma^2(b)}\frexp{-2\delta\frac{\mu(b)}{\sigma^2(b)}} \frac{\zeta(b)+\mu(b)}{2(r+\rho)\zeta(b)} \frexp{\frac{(z+\delta) (\mu(b)-\zeta(b))}{\sigma^2(b)}} \\
&\quad&+\frexp{-2\delta\frac{\mu(b)}{\sigma^2(b)}} \zeta(b)^{-1} \frexp{\frac{(z+\delta) (\mu(b)-\zeta(b))}{\sigma^2(b)}}\\
&=& \zeta(b)^{-1} \frexp{\frac{(z-\delta)\mu(b)-\abs{z-\delta} \zeta(b)}{\sigma^2(b)}} \\ 
&&\quad {}+\bigl[\zeta(b)^{-1}+\frac{\mu(b) (\zeta(b)+\mu(b))}{(r+\rho)\zeta(b)\sigma^2(b)} \bigr] \frexp{-2\delta\frac{\mu(b)}{\sigma^2(b)}} \frexp{\frac{(z+\delta) (\mu(b)-\zeta(b))}{\sigma^2(b)}} \;.
\end{eqnarray*}
Finally, we need the following result for developing an algorithm to solve equation \eqref{dynprorewritten}.
\begin{lemma}
    For $z \leq a$ it holds
    \begin{eqnarray*}
   \int_a^{\infty} w^0(\delta,z) \id \delta = \frac{\rho}{r+\rho}\frexp{-(a-z) \frac{r+\rho}{\theta-\eta}}
    \end{eqnarray*}
    and for $b>0$
    \begin{eqnarray*}
   &&\int_a^{\infty} w^b(\delta,z) \id \delta = \frac{\rho\sigma^2(b)}{\zeta(b)(\mu(b)+\zeta(b))} \frexp{-(a-z) \frac{\mu(b)+\zeta(b)}{\sigma^2(b)}}\\
   &\quad&+\rho \bigl[\frac{\sigma^2(b)}{\zeta(b)(\mu(b)+\zeta(b))}+\frac{\mu(b)}{(r+\rho)\zeta(b)}\bigr] \frexp{-(a-z) \frac{\mu(b)}{\sigma^2(b)}}\frexp{-(a+z) \frac{\zeta(b)}{\sigma^2(b)}}\;.
    \end{eqnarray*}
\hfill\qed
\end{lemma}

\subsubsection{Description of the Algorithm}

In order to calculate the value function numerically, we note that, due to the Riemann integrability of $v$ and $w^b$, equation \eqref{dynprorewritten} takes the form
\begin{eqnarray*}
v(z) &=& \onebr{z>d} + \inf_{b\in[0,1]} \bigl\{\int_0^\infty v(\delta)  w^b(\delta,z) \id \delta\bigr\} \\
&=& \onebr{z>d} + \inf_{b\in[0,1]}\bigl\{\sup_{\mathbf{p}=(p_0,\ldots,p_N)\in?P}\sum_{i=1}^{N} \frac{v(p_{i-1})+v(p_i)}{2} \int_{p_{i-1}}^{p_i} w^b(\delta,z) \id \delta + \int_a^\infty v(\delta)  w^b(\delta,z) \id \delta \bigr\}\;,
\end{eqnarray*}
where 
\[
?P = \{\mathbf{p}=(p_0,\ldots,p_N) : 0=p_0<p_1< \cdots<p_N=a, N\in!!N\} 
\]
is the set of partitions of the interval $[0,a]$ for $a>0$. Moreover, if $\tilde{?P}$ is the set of partitions of the interval $[0,1]$, we find
\begin{eqnarray*}
v(z) &=& \onebr{z>d} + \inf_{\mathbf{b}=(b_0,\ldots,b_M)\in\tilde{?P}}\min_{j=0,\ldots,M}\bigl\{\sup_{\mathbf{p}=(p_0,\ldots,p_N)\in?P}\sum_{i=1}^{N} \frac{v(p_{i-1})+v(p_i)}{2} \int_{p_{i-1}}^{p_i} w^{b_j}(\delta,z) \id \delta \\
&{}& \qquad + \int_a^\infty v(\delta)  w^{b_j}(\delta,z) \id \delta \bigr\}\;.
\end{eqnarray*}
We conclude for $k\in\{0,\ldots,n\}$:
\begin{eqnarray*}
v(z_k) &\approx& \onebr{z_k>d} + \min_{j=0,\ldots,M}\bigl\{\sum_{i=1}^{N} \frac{v(p_{i-1})+v(p_i)}{2} \int_{p_{i-1}}^{p_i} w^{b_j}(\delta,z_k) \id \delta + \frac{\rho+r}{r} \int_a^\infty w^{b_j}(\delta,z_k) \id \delta \bigr\}\;,
\end{eqnarray*}
if we choose the partitions $\mathbf{p},\mathbf{z}\in?P$ and $\mathbf{b}=(b_0,\ldots,b_M)\in\tilde{?P}$ sufficiently dense and $a>0$ large enough. We will now apply the idea from Remark \ref{rem:Banach}: We fix $N,M\in!!N$ and partitions $(z_0,\ldots,z_N)\in?P$ and $(b_0,\ldots,b_M)\in\tilde{?P}$. Choose $\mathbf{v_0}=(v_0^{(0)},\ldots,v_N^{(0)})\in[0,(\rho+r)/r]^{N+1}$ and define $\mathbf{v_{n+1}}=(v_0^{(n+1)},\ldots,v_N^{(n+1)})$ via
\[
v_k^{(n+1)} = \onebr{z_k>d} + \min_{j=0,\ldots,M}\bigl\{\sum_{i=1}^{N}\frac{v_{i-1}^{(n)}+v_i^{(n)}}{2} \int_{z_{i-1}}^{z_i} w^{b^j}(\delta,z_k) \id \delta + \frac{\rho+r}{r} \int_a^\infty  w^{b_j}(\delta,z_k) \id \delta \bigr\}\;.
\]
We stop the iteration as soon as we find $n\in!!N$ such that $\supnorm{\mathbf{v_{n+1}}-\mathbf{v_n}} <\epsilon$ for $\epsilon>0$ small.

\subsubsection{Numerical Results}

\begin{figure} 
\centering
\includegraphics[width=\textwidth]{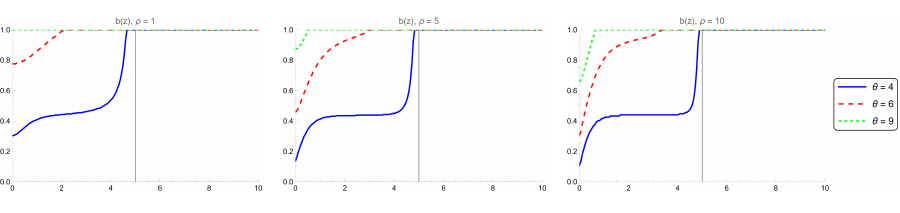}
\caption{Numerical results for the value function $v(z)$ and the corresponding optimal retention level $b(z)$ for different Poisson parameters $\rho\in\{1,5,10\}$.}
    \label{figure:numresults1}
\end{figure}

\begin{figure} 
\centering
\includegraphics[width=0.9\textwidth]{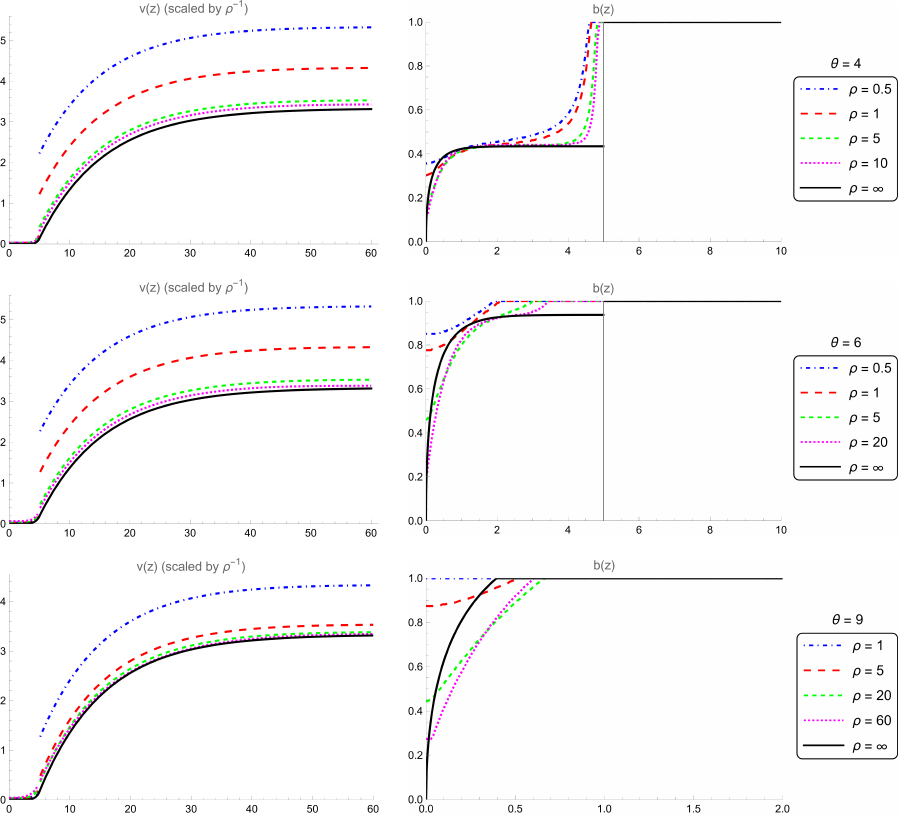}
\caption{Numerical results for the value function $v(z)$ and the corresponding optimal retention level $b(z)$ for different prices of reinsurance $\theta\in\{4,6,9\}$. The case ``$\rho = \infty$'' corresponds to the continuous time case of \cite{Brinker2022}.}
    \label{figure:numresults2}
\end{figure}

In the following, the results of the numerical analysis are illustrated and discussed. In Figure \ref{figure:numresults1} and \ref{figure:numresults2} we compare the value functions and corresponding retention levels for different prices of reinsurance $\theta$ and different Poisson parameters $\rho$. The selected parameter values are listed in Table \ref{table:parameters}. 

In general we can see, as observed in Lemma \ref{lemma:bound}, that the value
function is increasing and has a jump at the critical drawdown size $d$. In
Figure \ref{figure:numresults1} we perceive that the retention level increases
with the cost of reinsurance. Moreover, the retention level $b(z)$ increases as the
drawdown becomes larger until it constantly equals $1$. We observe that this point is reached sooner, when reinsurance is more expensive.
We can
compare these observations to the results of \cite{Brinker2022}, since the
continuous time model considered there can be interpreted as a limiting case
$\rho\to\infty$. In \cite{Brinker2022} the optimal strategy is given by a feedback strategy
$(\tilde{B}_t=\tilde{b}(\Delta_z^{\tilde{B}}(t)))_{t\geq 0}$, where the
minimiser $\tilde{b}(z)$ takes the form
\[
\tilde{b}(z) = \begin{cases}
\theta D\sigma^{-2}[1+?W(-\e^{-(1+zB/D})]\;, & z\in[0,z_0\wedge d]\;, \\
1\;, & z>z_0\wedge d\;.
\end{cases}
\]
Here, $?W(z)$ is the Lambert W function and $B,D>0$ and $z_0$ are
constants. In particular, it holds $\tilde{b}(z)=1$ for all $z>d$, which means
that we choose the strategy with maximal drift in the critical area in order to reenter the
non-critical area as quickly as possible. Comparing this with the results in Figures \ref{figure:numresults1} and \ref{figure:numresults2}, we can see that in our model there is also a certain drawdown size above which it is not optimal to buy
reinsurance, but this level is smaller than $d$ in the examples considered. Since we can only change the retention level at the observation times, i.e. in particular not at the moment when we enter the critical area, we will take countermeasures in advance when the drawdown approaches the critical drawdown level $d$.
The second difference we observe is the behaviour of the strategy for $z$
close to $ 0$: In \cite{Brinker2022} it holds $\tilde{b}(z) \to 0$ as $z\to
0$. Since the drift $\mu(b)=\eta-(1-b)\theta$ becomes negative for $b$ close to
$0$, this means that the maximum value is never exceeded. In our model here,
we can see that $\lim_{z\to 0} b^{*}(z)>0$ in all scenarios considered. However, in Figure \ref{figure:numresults2} we observe that this value decreases in $\rho$.
This behaviour can be explained in a similar way as above: In \cite{Brinker2022} the strategy can
be adjusted at any time, so that we can take immediate countermeasures if the
drawdown becomes too large. As mentioned before, in our model we have to wait for the
random times $\{T_k\}_{k\in!!N_0}$, such that a negative drift could lead to
excessive growth in the drawdown if the expected interarrival time $\rho^{-1}$ is large. All in all, we can see in Figure \ref{figure:numresults2} that the optimal retention level and the corresponding value function converge to the continuous time case of \cite{Brinker2022} as the expected interarrival time becomes smaller. Therefore the differences mentioned above disappear for $\rho\to\infty$.

\begin{table}[ht]
    \centering
    \begin{tabular}{c|c|c|c|c}
Parameter &  $\eta$ & $\sigma$ & $r$ & $d$ \\ \hline
Value & 3 & 2 & 0.3 & 5
    \end{tabular}
    \caption{Parameters of the insurer's safety loading, volatility, preference factor and critical drawdown size (from left to right).}
    \label{table:parameters}
\end{table}  


\subsection{Deterministic Observations}

In this case, the drawdown process is inspected at time points
$\{T_n\}_{n\in!!N_0}=(n T)_{n\in!!N_0}$ for a fixed interarrival time
$T>0$. This is the most intuitive case, where the insurer checks the size of
the drawdown at regular intervals (daily, weekly or quarterly, for example).

\subsubsection{Rewriting the Dynamic Programming Equation}

In view of Lemma \ref{lemma:bound}, the value function is bounded from above by $(1-\e^{-rT})^{-1}$. The dynamic programming equation \eqref{dynamicprogramming} can be easily rewritten as
\begin{eqnarray*}
    v(z) &=& \onebr{z>d} + \e^{-rT} \inf_{b\in[0,1]} \int_0^\infty v(\delta) \id F_{\Delta_z^b(T)}(\delta) = \onebr{z>d} + \e^{-rT} \inf_{b\in[0,1]} w^b(z)\;,
\end{eqnarray*}
where
\begin{equation*}
    w^{0}(z) = v(z+(\theta-\eta)T)
\end{equation*}
and
\begin{equation*}
    w^b(z) = \int_0^\infty v(\delta)  f_{\Delta_z^b(T)}(\delta) \id \delta
\end{equation*}
for $b>0$.
\subsubsection{Description of the Algorithm}
 We follow a similar approach as in the Poisson case. Let $\mathbf{z}=(z_0,\ldots,z_N)$ be a sufficiently dense partition of the interval $[0,a]$ for $a>0$ and $\mathbf{b}=(b_0,\ldots,b_M)$ one for the interval $[0,1]$. Then the approximated version of the dynamic programming equation becomes
 \begin{eqnarray*}
v(z_k) &\approx& \onebr{z_k>d} + \e^{-rT}  \min \bigl\{v(z_k+(\theta-\eta)T), \\
&\qquad& \min_{j=0,\ldots,M} \bigl\{\sum_{i=1}^N \frac{v(z_{i-1})+v(z_i)}{2}
(F_{\Delta_{z_{k}}^{b_j}(T)}(z_i)-F_{\Delta_{z_{k}}^{b_j}(T)}(z_{i-1})) + (1-\e^{-rT})^{-1}(1-F_{\Delta_{z_k}^{b_j}(T)}(a))\bigl\}\bigl\}\;,
\end{eqnarray*}
for $k\in\{0,\ldots,N\}$. Note that we need to make sure that $z_k+(\theta-\eta)  T$ is again a grid point or large enough that we can approximate the value function with the upper bound $v(z_k+(\theta-\eta)T)\approx(1-\e^{-rT})^{-1}$. We therefore construct the grid in the following way: Choose $n,K\in!!N$ and $\{z_i\}_{i=0}^n\subset[0,(\theta-\eta)T)$ with $0=z_0<z_1<\dots<z_n$. For $j\in\{1,\ldots,K\}$ and $i\in\{jn+1,\ldots,(j+1)n\}$ set
$
z_i = z_{i-jn}+(\theta-\eta)  T  j\;.
$
Thus, the sequence $\{z_i\}_{i=0}^{(K+1)n}$ forms a grid on the interval
$[0,(\theta-\eta)  T (K+1))$. For $K\in !!N$ large enough we can set $v(z_k)=(1-\e^{-rT})^{-1}$ for $k\in\{Kn+1,\ldots,(K+1)n\}$. Now one can numerically determine $v(z_k)$ for $k\in\{0,\ldots,Kn\}$ analogously to the Poisson case.

\subsubsection{Numerical Results}
\begin{figure}
        \centering
        \includegraphics[width=\textwidth]{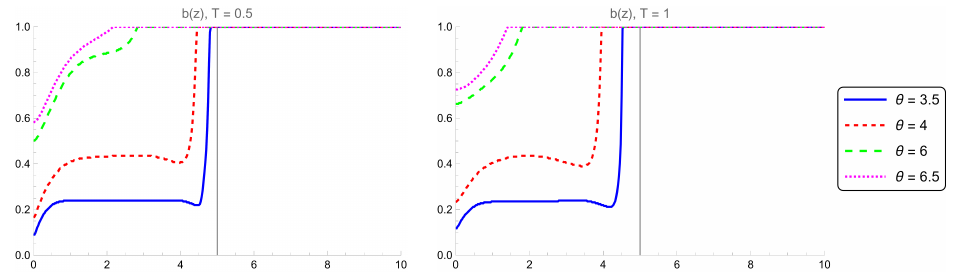}
        \caption{Numerical results for the value function $v(z)$ and the corresponding optimal retention level $b(z)$ for different interarrival times $T\in\{0.5,1\}$.}
        \label{figure:numresults3}
\end{figure}

\begin{figure}
        \centering
        \includegraphics[width=\textwidth]{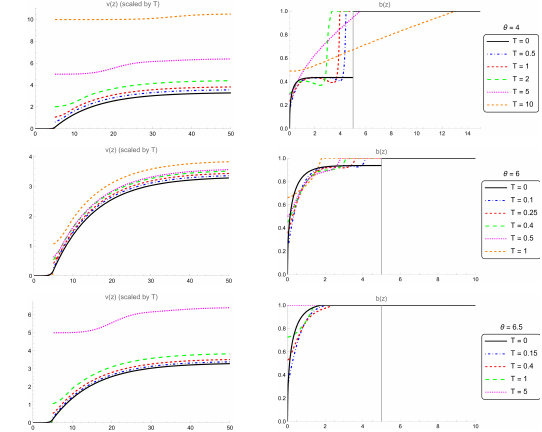}
                \caption{Numerical results for the value function $v(z)$ and the corresponding optimal retention level $b(z)$ for different prices of reinsurance $\theta\in\{4,6,6.5\}$. The case ``$T = 0$'' corresponds to the continuous time case of \cite{Brinker2022}.}
        \label{figure:numresults4}
\end{figure}

The results of the numerical study are illustrated in Figures
\ref{figure:numresults3} and \ref{figure:numresults4}. As in the Poisson case, we can see in Figure \ref{figure:numresults4} that the value function and the optimal strategy converge to the continuous time case of \cite{Brinker2022}. Moreover, the retention
level is again increasing in the cost of reinsurance $\theta$ (see Figure \ref{figure:numresults3}). However, it is not longer increasing in the drawdown size in every scenario considered and the behaviour of the strategy changes significantly for large interarrival times. In particular, buying reinsurance in the critical area is optimal in some cases. This can be explained as follows: In the model considered here, a critical drawdown is only penalized, when it is observed at the deterministic times $\{nT\}_{n\in!!N_0}$. If $T$ is large, there is more time to reenter the critical area, since we know the next observation time for sure. We can therefore choose a less risky strategy in order to increase the chance reaching the non-critical area before the next observation.

\newpage
\appendix

\section*{Appendix}

\begin{proof}[Proof of Lemma \ref{lemma:integralsinverseGaussian}]
We can rewrite the exponent as
\[
-st-\frac{(\alpha-\mu t)^2}{2t\sigma^2} = q -\frac{\tilde{\lambda}(t-\tilde{\mu})^2}{2\tilde{\mu}^2t}\;,
\]
where 
\[
\tilde{\lambda} = \frac{\alpha^2}{\sigma^2}\;, \qquad \tilde{\mu}=\frac{\abs{\alpha}}{\sqrt{\mu^2+2s\sigma^2}} \quad \text{and}\quad q= \frac{\alpha \mu-\abs{\alpha} \sqrt{\mu^2+2s\sigma^2}}{\sigma^2}\;.
\]
The function
\[
f(t) = \sqrt{\frac{\tilde{\lambda}}{2\pi t^3}}\frexp{-\frac{\tilde{\lambda}(t-\tilde{\mu})^2}{2\tilde{\mu}^2t}}\one_{t\geq 0}
\]
is the density function of the inverse Gaussian distribution $\text{IG}(\tilde{\mu},\tilde{\lambda})$ with parameters $\tilde{\mu}$ and $\tilde{\lambda}$. The expected value of a random variable $X\sim \text{IG}(\tilde{\mu},\tilde{\lambda})$ is given by $\E[X]=\tilde{\mu}$. 
We therefore conclude
\begin{eqnarray*}
\lefteqn{\int_0^\infty \e^{-st} \frac{1}{\sqrt{2\pi
t^3}\sigma}\frexp{-\frac{(\alpha-\mu t)^2}{2t\sigma^2}} \id t =
\frac{1}{\sigma\sqrt{\tilde{\lambda}}}\e^{q} \int_0^\infty \sqrt{\frac{\tilde{\lambda}}{2\pi
t^3}}\frexp{-\frac{\tilde{\lambda}(t-\tilde{\mu})^2}{2\tilde{\mu}^2t}} \id t}\hskip1cm && \hskip14cm
\null \\
&=& \frac{1}{\sigma\sqrt{\tilde{\lambda}}}\e^{q} = \abs{\alpha}^{-1} \frexp{\frac{\alpha \mu-\abs{\alpha} \sqrt{\mu^2+2s\sigma^2}}{\sigma^2}}
\end{eqnarray*}
and
\begin{eqnarray*}
\lefteqn{\int_0^\infty \e^{-st} \frac{1}{\sqrt{2\pi
t}\sigma}\frexp{-\frac{(\alpha-\mu t)^2}{2t\sigma^2}} \id t =
\frac{1}{\sigma\sqrt{\tilde{\lambda}}}\e^{q} \int_0^\infty t \sqrt{\frac{\tilde{\lambda}}{2\pi
t^3}}\frexp{-\frac{\tilde{\lambda}(t-\tilde{\mu})^2}{2\tilde{\mu}^2t}} \id t}\hskip1cm && \hskip14cm
\null \\
&=& \frac{1}{\sigma\sqrt{\tilde{\lambda}}}\e^{q} \tilde{\mu} = \frac{1}{\sqrt{\mu^2+2s\sigma^2}} \frexp{\frac{\alpha \mu-\abs{\alpha} \sqrt{\mu^2+2s\sigma^2}}{\sigma^2}}\;.
\end{eqnarray*}
Since
\begin{eqnarray*}
\lefteqn{\int_0^\infty \e^{-st} \frphi{\frac{\alpha-\mu t}{\sigma\sqrt{t}}} \id
t} \hskip1cm \\
&=& \frac{1}{s} \lim\limits_{t\to 0}\frphi{\frac{\alpha-\mu t}{\sigma\sqrt{t}}} -\frac{1}{s} \int_0^\infty \e^{-st} \biggl(\frac{\alpha}{2\sigma\sqrt{t^3}}+\frac{\mu}{2\sigma\sqrt{t}}\biggr) \frac{1}{\sqrt{2\pi}}\frexp{-\frac{(\alpha-\mu t)^2}{2t\sigma^2}} \id t \\
&=&\frac{1}{s} \biggl(\one_{\{\alpha>0\}}+\frac{1}{2}\one_{\{\alpha=0\}}\biggr)-\frac{\alpha}{2s} \int_0^\infty \e^{-st} \frac{1}{\sqrt{2\pi t^3}\sigma}\frexp{-\frac{(\alpha-\mu t)^2}{2t\sigma^2}} \id t \\
&\qquad&-\frac{\mu}{2s} \int_0^\infty \e^{-st} \frac{1}{\sqrt{2\pi t}\sigma}\frexp{-\frac{(\alpha-\mu t)^2}{2t\sigma^2}} \id t\;,
\end{eqnarray*}
the third identity follows from the calculations above.
\end{proof}
\bibliographystyle{abbrv}
\bibliography{Literatur}
\end{document}